\input amssym.def
\input amssym
\documentstyle[11pt]{article}
\begin{document}
\newtheorem{thm}{Theorem}[section]
\newtheorem{defin}[thm]{Definition}
\newtheorem{lemma}[thm]{Lemma}
\newtheorem{propo}[thm]{Proposition}
\newtheorem{cor}[thm]{Corollary}
\newtheorem{conj}[thm]{Conjecture}

\centerline{\huge \bf Braided Deformations of Monoidal}

\centerline{\huge \bf Categories and Vassiliev Invariants}
\bigskip

\centerline{\parbox{50mm}{David N. Yetter \\ Department of
Mathematics \\ Kansas State University \\ Manhattan, KS 66506}\footnote{
Supported by NSF Grant \# DMS-9504423.  \\ This paper has been submitted
to the proceedings of the Workshop on Higher Dimensional Category Theory
and Mathematical Physics, Northwestern University, March 27-28, 1997}}
\medskip

\section{Introduction}

It is the purpose of this paper to provide an exposition of a cohomological
deformation theory for braided monoidal categories,
an exposition and proof of a very
general theorem of the form ``all quantum invariants are Vassiliev
invariants,'' and to relate
to the existence of the Vassiliev invariants obtained via the theorem to
the vanishing of obstructions in the deformation theory.

Our deformation theory for braided monoidal categories will be based
on a deformation theory for monoidal functors.  Remarkably, the deformation
theory of monoidal functors, unlike the deformation theory of monoidal
categories developed in \cite{CYdef} shares many of
the properties of the Gerstenhaber\cite{G1,G2,GS}
deformation theory of algebras, though the relationship between the two 
theories is unclear.\footnote{It has been recently brought to the
author's attention that some aspects of this work have been discovered
independently by Davydov \cite{Dav}. Davydov's work does not however deal
with ``higher-order infinitesimal deformations'', and thus the full
parellel with Gerstenhaber's deformation theory is unrealized.}
In this present paper, the deformation theory of
monoidal functors is worked out only in the case where the monoidal
structure of the source and target categories are undeformed.

 Previous theorems of the ``all quantum invariants are Vassiliev''
type have involved
definitions of
``quantum invariants'' either in terms of the well-known knot polynomial
of Jones \cite{Jones}, HOMFLY \cite{HOMFLY}, or Kauffman \cite{Kauffman},
as in the case of Birman and Lin \cite{BL} and Stanford \cite{Stanford}, 
or more generally in terms
of constructions tied explicitly to a (simple) Lie algebra, as in the case of
Piunikhin \cite{Piu}.
Still, the notion of ``quantum invariants'' used here is not the most
general conceivable, so some preliminaries are in order.  The most
general reasonable notion of ``quantum invariant'' of (framed)
knots and links is
what will here be called a ``functorial invariant'', that is, an $R$-valued
invariant ($R$ some commutative ring) which arises by choosing an object  
in an $R$-linear tortile (a.k.a. ribbon) category (with $End({\bf I}) = R$,
where {\bf I} is the monoidal identity object), and considering the
image of all framed links under the functor induced by the coherence theorem
of Shum \cite{Shum}.  If we wish to deal with unframed links in this 
context, we must find a normalization procedure. (One will always be
available if the chosen object is simple.)

It is not the case that all functorial invariants are Vassiliev invariants:
the ``counting invariants'' associated to finite groups are functorial
(cf. Freyd and Yetter \cite{FY.BCCC} and 
Yetter \cite{Y.FG}), but as shown by Altschuler \cite{Alt} are not
generally Vassiliev.   Indeed, even the functorial invariants usually
shown to ``be'' Vassiliev invariants are not themselves Vassiliev invariants,
but rather become a sequence of Vassiliev invariants when a suitable
change of variable by substitution of a power series is applied. Our theorem
will be stated in terms of those functorial invariants which arise from
tortile (ribbon) deformations of a rigid symmetric monoidal category.

Our setting is more general than that used for previous 
results of this type,
and, we have reason to conjecture, may be the most general setting in which
a theorem of this type can be proved. We will point to examples of knot
invariants already existing in the literature to show
the greater generality. More importantly, by proving the 
result in a purely categorical setting, we are able to relate the Vassiliev
invariants thus obtained to a new cohomology theory which both classifies
and provides the obstruction theory for deformations of braidings on
tensor categories. Once the relationship between our deformation theory
and the Gerstenhaber deformation theory of algebras is clearer, it is
reasonable to think that we will gain an improved understanding of the
relation between quantum groups and Vassiliev theory.

 We begin with a discussion of the
categorical machinery necessary for the statement and proof of the theorem.

\section{Tortile Categories and Deformations of Tensor Categories and Functors}

The setting in which our theorem is proved involves a class of categories 
with structure which includes both the categories of representations of
quantum groups and various categories of tangles.  To make this explicit:

\begin{defin}
A {\em monoidal category} $\cal C$ is a category $\cal C$ equipped
with a functor $\otimes : {\cal C} \times {\cal C} \rightarrow {\cal C}$
and an object $I$, together with natural isomorphisms 
$\alpha : \otimes (\otimes \times 1_{\cal C}) \Rightarrow \otimes (1_{\cal C}
\times \otimes)$, $\rho : \otimes I \Rightarrow 1_{\cal C}$, and 
$\lambda : I\otimes \Rightarrow 1_{\cal C}$, satisfying the
pentagon, triangle, and bigon ($\rho_I = \lambda_I$) coherence conditions
(cf. \cite{CWM}).  A {\em tensor category over $K$}, for $K$ some field
(or commutative ring) is a monoidal abelian category linear over $K$, with
all $-\otimes X$ and $X\otimes -$ exact.  Similarly, a {\em semigroupal
category} is a category equipped with only $\otimes$ and $\alpha$ satisfying
the usual pentagon.
\end{defin}

\begin{defin}
A {\em (strong) monoidal functor} $F:{\cal C}\rightarrow {\cal D}$ between
two monoidal categories $\cal C$ and $\cal D$ is a functor $F$ between
the underlying categories, equipped with a natural isomorphism

\[ \tilde{F}: F(-\otimes -)\rightarrow F(-)\otimes F(-) \]

\noindent and an isomorphism
$F_0:F(I)\rightarrow I$, satisfying the hexagon and two squares
\end{defin}

\begin{figure}[h]

\setlength{\unitlength}{0.00083300in}%
\begin{picture}(6300,5055)(301,-4372)
\thicklines
\put(3150,389){\vector( 1, 0){1500}}
\put(2326,164){\vector( 0,-1){900}}
\put(2326,-1261){\vector( 0,-1){750}}
\put(5251,164){\vector( 0,-1){825}}
\put(5251,-1186){\vector( 0,-1){825}}
\put(3450,-2236){\vector( 1, 0){1025}}
\put(1801,314){\makebox(0,0)[lb]{\raisebox
{0pt}[0pt][0pt]{\twlrm $F((A\otimes  B)\otimes  C)$}}}
\put(1801,-1111){\makebox(0,0)[lb]{\raisebox
{0pt}[0pt][0pt]{\twlrm $F(A\otimes  B)\otimes  F(C)$}}}
\put(3526,539){\makebox(0,0)[lb]{\raisebox{0pt}[0pt][0pt]{\twlrm $F(\alpha)$}}}
\put(1576,-361){\makebox(0,0)[lb]{\raisebox
{0pt}[0pt][0pt]{\twlrm $\tilde{F}_{A\otimes  B,C}$}}}
\put(5476,-436){\makebox(0,0)[lb]{\raisebox
{0pt}[0pt][0pt]{\twlrm $\tilde{F}_{A,B\otimes  C}$}}}
\put(1125,-1711){\makebox(0,0)[lb]{\raisebox
{0pt}[0pt][0pt]{\twlrm $\tilde{F}_{A,B}\otimes  F(C)$}}}
\put(4876,314){\makebox(0,0)[lb]{\raisebox
{0pt}[0pt][0pt]{\twlrm $F(A\otimes  (B\otimes  C))$}}}
\put(4801,-1036){\makebox(0,0)[lb]{\raisebox
{0pt}[0pt][0pt]{\twlrm $F(A)\otimes  F(B\otimes  C)$}}}
\put(4726,-2311){\makebox(0,0)[lb]{\raisebox
{0pt}[0pt][0pt]{\twlrm $F(A)\otimes  (F(B)\otimes  F(C))$}}}
\put(1576,-2311){\makebox(0,0)[lb]{\raisebox
{0pt}[0pt][0pt]{\twlrm $(F(A)\otimes  F(B))\otimes  F(C)$}}}
\put(3751,-2536){\makebox(0,0)[lb]{\raisebox{0pt}[0pt][0pt]{\twlrm $\alpha$}}}
\put(5476,-1636){\makebox(0,0)[lb]{\raisebox
{0pt}[0pt][0pt]{\twlrm $F(A)\otimes \tilde{F}_{B,C}$}}}
\put(976,-3361){\vector( 0,-1){450}}
\put(1351,-3136){\line( 0, 1){  0}}
\put(1575,-3136){\vector( 1, 0){900}}
\put(2776,-3811){\vector( 0, 1){450}}
\put(1725,-4036){\vector( 1, 0){525}}
\put(676,-3211){\makebox(0,0)[lb]{\raisebox
{0pt}[0pt][0pt]{\twlrm $F(I\otimes A)$}}}
\put(2626,-3211){\makebox(0,0)[lb]{\raisebox{0pt}[0pt][0pt]{\twlrm $F(A)$}}}
\put(2551,-4111){\makebox(0,0)[lb]{\raisebox
{0pt}[0pt][0pt]{\twlrm $I\otimes  F(A)$}}}
\put(601,-4111){\makebox(0,0)[lb]{\raisebox
{0pt}[0pt][0pt]{\twlrm $F(I)\otimes  F(A)$}}}
\put(1726,-2986){\makebox(0,0)[lb]{\raisebox
{0pt}[0pt][0pt]{\twlrm $F(\lambda)$}}}
\put(2926,-3736){\makebox(0,0)[lb]{\raisebox{0pt}[0pt][0pt]{\twlrm $\lambda$}}}
\put(1501,-4336){\makebox(0,0)[lb]{\raisebox
{0pt}[0pt][0pt]{\twlrm $F_0\otimes  F(A)$}}}
\put(301,-3661){\makebox(0,0)[lb]{\raisebox
{0pt}[0pt][0pt]{\twlrm $\tilde{F}_{I,A}$}}}
\put(4651,-3361){\vector( 0,-1){450}}
\put(5026,-3136){\line( 0, 1){  0}}
\put(5250,-3136){\vector( 1, 0){900}}
\put(6451,-3811){\vector( 0, 1){450}}
\put(5400,-4036){\vector( 1, 0){525}}
\put(4351,-3211){\makebox(0,0)[lb]{\raisebox
{0pt}[0pt][0pt]{\twlrm $F(A\otimes I)$}}}
\put(6301,-3211){\makebox(0,0)[lb]{\raisebox{0pt}[0pt][0pt]{\twlrm $F(A)$}}}
\put(6226,-4111){\makebox(0,0)[lb]{\raisebox
{0pt}[0pt][0pt]{\twlrm $F(A)\otimes I$}}}
\put(4276,-4111){\makebox(0,0)[lb]{\raisebox
{0pt}[0pt][0pt]{\twlrm $F(A)\otimes F(I)$}}}
\put(5401,-2986){\makebox(0,0)[lb]{\raisebox{0pt}[0pt][0pt]{\twlrm $F(\rho)$}}}
\put(6601,-3736){\makebox(0,0)[lb]{\raisebox{0pt}[0pt][0pt]{\twlrm $\rho$}}}
\put(5176,-4336){\makebox(0,0)[lb]{\raisebox
{0pt}[0pt][0pt]{\twlrm $F(A)\otimes  F_0$}}}
\put(3976,-3661){\makebox(0,0)[lb]{\raisebox
{0pt}[0pt][0pt]{\twlrm $\tilde{F}_{A,I}$}}}
\end{picture}

\caption{Coherence conditions for a monoidal functor}
\end{figure}

\clearpage
\noindent{\em A pair of a functor $F$ and a natural isomorphism 
$\tilde{F}$ (without the isomorphism $F_0$
as above will be called a {\em semigroupal functor}.}
\smallskip

\begin{defin}
A {\em monoidal natural transformation} is a natural transformation
$\phi:F\Rightarrow G$ between monoidal functors which satisfies

\[ \tilde{G}_{A,B}(\phi_{A\otimes B}) = \phi_A \otimes \phi_B (\tilde{F}_{A,B})
\]

\noindent and $F_0 = G_0(\phi_I)$.  A {\em semigroupal natural transformation}
between semigroupal functors is defined similarly.
\end{defin}

We will be concerned with monoidal categories with additional structure.

\begin{defin}
A {\em braided monoidal category} is a monoidal category equipped with a
monoidal natural isomorphism $\sigma: \otimes \Rightarrow \otimes(tw)$,
where $tw:  {\cal C} \times {\cal C} \rightarrow {\cal C} \times {\cal C}$
is the ``twist functor'' ($tw(f,g) = (g,f)$).  And satisfying
\end{defin}

\begin{figure}[h]
\setlength{\unitlength}{0.0125in}%
\begin{center} \begin{picture}(335,161)(50,615)
\thinlines
\put( 90,710){\vector( 4, 3){ 40}}
\put(205,760){\vector( 1, 0){ 75}}
\put(205,640){\vector( 1, 0){ 75}}
\put( 90,685){\vector( 1,-1){ 35}}
\put(345,750){\vector( 1,-1){ 40}}
\put(340,650){\vector( 4, 3){ 40}}
\put( 50,690){\makebox(0,0)[lb]{\raisebox{0pt}[0pt][0pt]
{\twlrm $(A \otimes B) \otimes C$}}}
\put(125,755){\makebox(0,0)[lb]{\raisebox{0pt}[0pt][0pt]
{\twlrm $A \otimes (B \otimes C)$}}}
\put(295,755){\makebox(0,0)[lb]{\raisebox{0pt}[0pt][0pt]
{\twlrm $(B \otimes C) \otimes A$}}}
\put(130,635){\makebox(0,0)[lb]{\raisebox{0pt}[0pt][0pt]
{\twlrm $(B \otimes A) \otimes C$}}}
\put(290,635){\makebox(0,0)[lb]{\raisebox{0pt}[0pt][0pt]
{\twlrm $B \otimes (A \otimes C)$}}}
\put(365,690){\makebox(0,0)[lb]{\raisebox{0pt}[0pt][0pt]
{\twlrm $B \otimes (C \otimes A)$}}}
\put( 85,730){\makebox(0,0)[lb]{\raisebox{0pt}[0pt][0pt]
{\twlrm $\alpha$}}}
\put(235,770){\makebox(0,0)[lb]{\raisebox{0pt}[0pt][0pt]
{\twlrm $\sigma^{\pm 1}$}}}
\put(380,735){\makebox(0,0)[lb]{\raisebox{0pt}[0pt][0pt]
{\twlrm $\alpha$}}}
\put( 60,655){\makebox(0,0)[lb]{\raisebox{0pt}[0pt][0pt]
{\twlrm $\sigma^{\pm 1} \otimes C$}}}
\put(230,615){\makebox(0,0)[lb]{\raisebox{0pt}[0pt][0pt]
 {\twlrm $\alpha$}}}
\put(380,655){\makebox(0,0)[lb]{\raisebox{0pt}[0pt][0pt]
{\twlrm $B \otimes \sigma^{\pm 1}$}}}
\end{picture} \end{center}

\caption{The hexagon}
\end{figure}

\noindent{\em a braided monoidal category is a 
{\em symmetric monoidal category}
if the components of $\sigma$ satisfy 
$\sigma_{B,A} (\sigma_{A,B}) = 1_{A\otimes B}$ for all objects $A$ and $B$.
}

\begin{defin}
A {\em right dual} to an object $X$ in a monoidal category $\cal C$ is
an object $X^\ast$ equipped with maps 
$\epsilon: X \otimes X^\ast \rightarrow I$ and
$\eta: I \rightarrow X^\ast \otimes X$ such that the compositions

\[ X \stackrel{\rho^{-1}}{\rightarrow} X \otimes I 
\stackrel{X\otimes \eta}{\rightarrow} 
X \otimes (X^\ast \otimes X) \stackrel{\alpha^{-1}}{\rightarrow} 
(X\otimes X^\ast )\otimes
X \stackrel{\epsilon \otimes X}{\rightarrow} I\otimes X 
\stackrel{\lambda}{\rightarrow} X \]

\noindent and

\[ X^\ast \stackrel{\lambda^{-1}}{\rightarrow} I \otimes X^\ast 
\stackrel{\eta \otimes X^\ast }{\rightarrow} 
(X^\ast \otimes X) \otimes X^\ast \stackrel{\alpha}{\rightarrow} 
X^\ast \otimes 
(X \otimes X^\ast) \stackrel{X^\ast \otimes \epsilon}{\rightarrow} 
X^\ast \otimes I 
\stackrel{\rho}{\rightarrow} X^\ast \]

\noindent are both identity maps.

A {\em left dual} $^\ast X$ is defined similarly, 
but with the objects placed on opposite
sides of the monoidal product.
\end{defin}

This type of duality is an abstraction from the sort of duality which exists
in categories of finite dimensional vectorspaces. It is not hard to show
that the canonical isomorphism from the second dual of a vectorspace
to the space generalizes to give canonical isomorphisms 
$k: \ast(X^\ast) \rightarrow X$ and $\kappa: (^\ast X)^\ast \rightarrow X$.
In general, however, there are not necessarily any maps from
$X^{\ast \ast}$ or $^{\ast \ast} X$ to X (cf. \cite{FY.cohere}).
In cases where every object admits a right (resp. left) dual, it is easy to
show that a choice of right (resp. left) dual for every object extends to
a contravariant functor, whose application to maps will be denoted
$f^\ast$ (resp. $^\ast f$), and that the canonical maps noted above
become natural isomorphisms between the compositions of these functors
and the identity functor.  Likewise, it is easy to show that
$(A\otimes B)^\ast$ is canonically isomorphic to $B^\ast \otimes A^\ast$, 
and similarly for left duals.

In the case of a braided monoidal category every right dual
is also a left dual, in general in a non-canonical way
(cf. \cite{FY.cohere}). In symmetric monoidal
categories, we return to the familiar: right duals are canonically
left duals. For non-symmetric braided monoidal categories, the structure
will be canonical only in the presence of additional structure.  As we will be 
concerned only with the case of categories in which all objects admit
duals, we make

\begin{defin}
A braided monoidal category $\cal C$ is {\em tortile (or ribbon)} 
if all objects
admit right duals, and it is equipped with a natural transformation
$\theta: 1_{\cal C} \Rightarrow 1_{\cal C}$, which satisfies

\[ \theta_{A\otimes B} = \sigma_{B,A}(\sigma_{A,B}(\theta_A \otimes \theta_B))
\]

\noindent and

\[ \theta_{A^\ast} = \theta_A^\ast \] 
\end{defin}

\begin{defin}
A symmetric monoidal category $\cal C$ is {\em rigid} if
all objects admit (right) duals.
\end{defin}

For details of the canonicity of the correspondence between right
and left duals in the tortile case see \cite{Y.FTTD}.
 
The importance of tortile categories for knot theory is provided by the
extremely beautiful coherence theorem of Shum \cite{Shum}:

\begin{thm} The tortile category freely generated by a single object
$X$ is monoidally equivalent to the category of framed tangles $\cal FT$.
\end{thm}

To be precise $\cal FT$ is the category whose objects are framed, oriented
(i.e. signed) finite sets of points in $(0,1)^2$, and whose arrows
are framed oriented compact
1-submanifolds of $[0,1]^3$ whose intersection with the boundary is
normal and consists of finite sets of points
lying in $(0,1)^2\times \{0\}$ (the source)
and $(0,1)^2\times \{1\}$ (the target), modulo ambient isotopy rel boundary.
Composition and $\otimes$ are defined by pasting on the third and second
coordinates respectively and rescaling.  Associativity and unit 
transformations are given by the obvious ascending submanifolds, while
the braiding is given by ascending submanifolds which pass one in front
of the other, and duality is given on objects by mirror-imaging in the second
coordinate, and reversing all signs, with half curves-of-rotation as 
structure maps (cf. \cite{Y.FTTD}).

This theorem gives rise to the following notion:

\begin{defin}
A {\em functorial invariant of framed links} $P_X$ is an $End(I)$-valued
invariant of framed links, where $I$ is the unit object in a tortile
category $\cal C$, obtained by choosing an object $X$ in $\cal C$,
and considering the restriction of the functor
 
\[  {\cal FT} \stackrel{G}{\longrightarrow} {\cal F} 
\stackrel{F_X}{\longrightarrow}
{\cal C}\]

\noindent where $\cal F$ is the free tortile category on one object
generator, $G$ is the equivalence functor from Shum's Coherence Theorem, and
$F_X$ is the functor induced by freeness which maps the generator to $X$.
\end{defin}

Numerous link invariants discovered since the mid 1980's
are functorial invariants, including the Jones polynomial, and Zariski
dense sets of values for the HOMFLY and Kauffman polynomials (the values
``arising from quantum groups''), as is Conway's normalization of the
classical Alexander polynomial.  

The standard algebraic construction of these functorial invariants begins
with the deformation of a universal enveloping algebra for a (simple) 
Lie algebra to produce a ``quantum group'', followed by the construction
of the (tortile) category of its representations. 

In \cite{Y.FTTD} the
author gave an exposition of a theorem of Deligne which dealt with 
direct deformation of a tensor category (starting with a Tannakian
category).  We will be concerned with the deformation theory of tensor 
categories and functors and its application to Vassiliev theory, so we review 
the necessary definitions:

\begin{defin}
  An {\em infinitesimal deformation} 
of a $K$-linear tensor category $\cal C$ over
an Artinian local $K$-algebra
$R$ is an $R$-linear tensor category $\tilde{\cal C}$
with the same objects as $\cal C$, but with $Hom_{\tilde{\cal C}} (a,b) =
Hom_{\cal C} (a,b)\otimes_K R$, and composition extended by bilinearity, 
and for which the structure map(s) $\alpha$ ($\rho$, and $\lambda$, and
if applicable $\sigma$, $\eta$ and $\epsilon$)
reduce mod $\frak m$ to the structure maps for $\cal C$, where 
$\frak m$ is the maximal ideal of $R$.  A deformation over 
$K[\epsilon]/<\epsilon ^{n+1}>$ is an {\em $n^{th}$ order deformation}.

Similarly an {\em $\frak m$-adic deformation} of $\cal C$ over an
$\frak m$-adically complete local $K$-algebra $R$ is an $R$-linear 
tensor category $\tilde{\cal C}$
with the same objects as $\cal C$, but with $Hom_{\tilde{\cal C}} (a,b) =
Hom_{\cal C} (a,b)\widehat{\otimes_K} R$, and composition extended 
by bilinearity and
continuity, 
and for which the structure map(s) $\alpha$  ($\rho$, and $\lambda$,
and if applicable $\sigma$, $\eta$ and $\epsilon$)
reduce mod $\frak m$ to the structure maps for $\cal C$, where 
$\frak m$ is the maximal ideal of $R$. (Here $\widehat{\otimes_K}$ is
the $\frak m$-adic completion of the ordinary tensor product.)  
An $\frak m$-adic deformation over 
$K[[x]]$ is a{\em formal series deformation}.

Two deformations  (in any of the above senses) are {\em equivalent}
if there exists a monoidal or semigroupal functor, whose underlying functor is
the identity, and whose structure maps reduce mod $\frak m$ to identity
maps from one to the other. {\em The trivial deformation} of $\cal C$ is the
deformation whose structure maps are those of $\cal C$. More generally,
a deformation is {\em trivial} if it is equivalent to a trivial deformation.

Finally, if $K = \Bbb C$ (or $\Bbb R$), and all hom-spaces in
$\cal C$ are finite dimensional, a {\em finite deformation} of $\cal C$ is
a $K$-linear tensor category with the same objects and maps 
as $\cal C$, but with
structure maps given by the structure maps of a formal series deformation
evaluated at $x = \xi$ for some $\xi \in K$ such that the formal series
defining all of the structure maps converge at $\xi$.

\end{defin}

Similarly, we can define deformations of monoidal and semigroupal
functors targetted at tensor categories by

\begin{defin}
Given a  $K$-linear tensor category $\cal D$.
  An {\em infinitesimal deformation} 
of a monoidal (resp. semigroupal) functor $F:{\cal C}\rightarrow {\cal D}$ 
 over
an Artinian local $K$-algebra
$R$ is a monoidal (resp. semigroupal)
functor $F^\prime:{\cal C}\rightarrow \tilde{\cal D}$
where the $R$-linear tensor category $\tilde{\cal D}$ 
is an infinitesimial deformation of $\cal D$ over $R$ 
and for which the structure maps $\tilde{F^\prime}$ and $F_0^\prime$ (as well
as those for $\cal D$)
reduce mod $\frak m$ to the structure maps for $F$ (and $\cal D$), where 
$\frak m$ is the maximal ideal of $R$.  A deformation over 
$K[\epsilon]/<\epsilon ^{n+1}>$ is an {\em $n^{th}$ order deformation}.

Similarly an {\em $\frak m$-adic deformation} of $F$ over an
$\frak m$-adically complete local $K$-algebra $R$ is a monoidal or semigroupal
functor targetted at the $R$-linear 
tensor category $\tilde{\cal D}$ is an $\frak m$-adic deformation of $\cal D$
and for which the structure maps 
reduce mod $\frak m$ to the structure maps for $F$ and $\cal D$, where 
$\frak m$ is the maximal ideal of $R$. (Here $\widehat{\otimes_K}$ is
the $\frak m$-adic completion of the ordinary tensor product.)  
An $\frak m$-adic deformation over 
$K[[x]]$ is a {\em formal series deformation}.

Two deformations  (in any of the above senses) are {\em equivalent}
if there exists a monoidal or semigroupal 
natural isomorphism between them which
reduces mod $\frak m$ to identity
maps. {\em The trivial deformation} of $F$ is the deformation whose
structure maps are those of $F$.  More generally, a deformation is 
{\em trivial} if it is equivalent to the trivial deformation.

A deformation in any of the above senses will be called {\em purely functorial}
whenever the target category $\tilde{\cal D}$ is the trivial deformation
of $\cal D$ .

\end{defin}

In what follows, the category or functor  we deform will generally have a more
restrictive structure than its deformations. Thus,
we will speak of a braided deformation or a tortile deformation of a
rigid symmetric category or a semigroupal deformation of a monoidal 
functor, meaning
one in which only the specified coherence conditions are satisfied.

Deligne's theorem \cite{Y.FTTD} then says

\begin{thm}
Any braided deformation of a Tannakian category $\cal C$ admits a unique
tortile structure, and is thus a tortile deformation.
\end{thm}
 
In fact the proof given in \cite{Y.FTTD} of this result will carry a
more general result once it is observed that the use of intertwining
matrices is unnecessary, and that it suffices to note that the map in
question is the sum of the identity map with $\frak m$-multiples of maps
in $\cal C$.

\begin{thm}
Any braided monoidal
deformation of a rigid symmetric $K$-linear tensor category 
$\cal C$
for $K$ any field (char $k \neq 2$) admits a unique tortile structure, and
is thus a tortile deformation.
\end{thm}

These theorems are quite helpful, since they show that even for topological
applications where the full tortile structure is required, it is sufficient
to understand the structure of braided deformations.

We will also have cause to consider the notion of a reduction of a deformation.

\begin{defin}
If $\cal C$ is a $K$-linear semigroupal category, $R$ an Artinian local ring
(resp. an $\frak m$-adically complete local ring), $\tilde{\cal C}$ 
a deformation ($\frak m$-adic deformation) of $\cal C$ over $R$, and
$J$ an ideal of $R$, then the {\em reduction $\bmod J$ of $\tilde{\cal C}$},
denoted $\tilde{\cal C}/J$, is the category whose hom-spaces
are $Hom_{\tilde{\cal C}}(X,Y)\otimes_R R/J$ with the obvious composition
induced by that on $\tilde{\cal C}$, and structure maps obtained as the
tensor product of the structure maps of $\tilde{\cal C}$ with $1 \in R/J$.
\end{defin}
 
\begin{defin}
If $\cal C$ is a $K$-linear semigroupal category, $F:{\cal D}\rightarrow
{\cal C}$  a semigroupal functor, $R$ an Artinian local ring
(resp. an $\frak m$-adically complete local ring), 
$F^\prime:{\cal D}\rightarrow
\tilde{\cal C}$ 
a deformation (resp. $\frak m$-adic deformation) of $F$ over $R$, and
$J$ an ideal of $R$, then the {\em reduction $\bmod J$ of $F^\prime$},
denoted $F^\prime /J$, is the monoidal functor targetted at $\tilde{\cal C}/J$
 and structure maps obtained as the
tensor product of the structure maps of $\tilde{\cal C}$ with $1 \in R/J$.
\end{defin} 

We will then need the following:

\begin{lemma} \label{reduction}
If $\tilde{\cal C}$ is an $n^{th}$ order deformation of a $K$-linear
semigroupal category $\cal C$ (resp. a formal
series deformation) (perhaps with additional
structure, e.g., monoidal, braided or tortile), then for $k < n$ (resp. for any
$k \in \Bbb N$), $\tilde{\cal C}/<\epsilon^{k+1}>$ is
a $k^{th}$ order deformation of $\cal C$ with the same additional
structure.  Similarly if $F^\prime$ is an $n^{th}$ order deformation of
a semigroupal 
functor targetted at $\cal C$ (resp. a formal series deformation),
then for any $k < n$ (resp. any $k \in {\Bbb N}$) $F^\prime/<\epsilon^{k+1}>$
is a $k^{th}$ order deformation of $F$.  If, moreover, $F^\prime$ is 
purely functorial, so is $F^\prime/<\epsilon^{k+1}>$.
\end{lemma}

\noindent {\bf proof:} It suffices to observe that the hom-spaces
of the reduction are those of a $k^{th}$ order deformation, and that
coherence conditions which hold before reduction $\bmod \; J$ 
will still hold after
reduction.
\smallskip

Until recently the notion of a deformation of a tensor category was
a definition without any attendant theorems, examples coming primarily
from categories of representations of a deformed Hopf algebra (quantum
group).  In \cite{CYdef} an exposition was given of the deformation of
tensor categories without braiding in terms of a cohomology theory for
tensor categories.  This theory is somewhat unsatisfactory in that it
is unclear whether the obstructions are of a purely cohomological nature.

For our purposes, however, it is much more interesting to consider the
deformation theory of tensor functors rather than tensor categories.
On the one hand, the obstruction theory will prove to be purely
cohomological, and on the other, it provides the correct notion of
deformation for braidings thanks to the following results of Joyal and
Street \cite{JS.BTC}:

\begin{defin}
A multiplication on a monoidal category $\cal C$ is a (strong) monoidal functor
$(\Phi:{\cal C}\times {\cal C}\rightarrow {\cal C}, \tilde{\Phi}, \Phi_0)$ 
(usually denoted $\Phi$ by abuse of notation)
together with monoidal natural
isomorphisms  ${\frak r}:\Phi(Id_{\cal C}, I) \Rightarrow
Id_{\cal C}$ and ${\frak l}:\Phi(I, Id_{\cal C}) \Rightarrow
Id_{\cal C}$
\end{defin}

\begin{thm} \label{mult.from.braiding}
In a monoidal category $\cal C$, a family of arrows

\[ \sigma_{A,B}:A\otimes B \longrightarrow B\otimes A \]

\noindent
is a braiding if and only if the following define a multiplication $\Phi$
on $\cal C$:  $\Phi = \otimes$, $\Phi_0 = \rho_I^{-1}$, ${\frak r} = \rho$,
${\frak l} = \lambda$ and 

\[ \tilde{\Phi}_{(A,A^\prime),(B,B^\prime)}
  = \lceil (1\otimes \sigma)\otimes 1
\rceil :(A\otimes A^\prime)\otimes (B\otimes B^\prime)\longrightarrow
(A\otimes B)\otimes (A^\prime \otimes B^\prime) \]

\end{thm}

Here, and in what follows, the notation $\lceil \;\; \rceil$ applied
to a map in a monoidal category indicates pre- and post- composition by the
unique map given by Mac Lane's coherence theorem \cite{CWM}
 as a composite of 
``prolongations'' of the structural transformations $\alpha$, $\rho$, and
$\lambda$ so that the source and targets will be as specified.  We will
also apply this notation to sequences of maps, in which case, Mac Lane
coherence maps may be inserted between each pair to provide a well-defined
composition, and to both single maps and sequences of maps in settings 
where a monoidal functor is involved.  In this
latter setting the
coherence map(s) may also involve structural transformations from the 
monoidal functor (see Epstein \cite{epstein} for the theorems which justify
use of this notation in the case of monoidal functors).

A converse theorem also due to Joyal and Street states

\begin{thm}
For any multiplication $\Phi$
on a monoidal category $\cal C$, a braiding $\sigma$ for $\cal C$ is defined 
by

\[ \sigma_{A,B} = {\frak l}\otimes {\frak r}(\tilde{\Phi}^{-1}(\Phi(\rho^{-1},
\lambda^{-1})(\Phi(\lambda,\rho)(\tilde{\Phi}
({\frak l}^{-1}\otimes 
{\frak r}^{-1}))))) \]

\noindent
The multiplication obtained from $\sigma$ via the previous theorem is
isomorphic (an the obvious sense) to $\Phi$, and if $\tau$ is any braiding
on $\cal C$, and $\Phi$ is the multiplication obtained from $\tau$ via
the previous theorem, then $\sigma = \tau$.
\end{thm}

As noted above, the 
importance of these theorems for us lies in the fact that they reduce
the question of studying braidings to the study of monoidal functors, and
it is from this point of view that we will examine the deformation theory
of braided tensor categories.

\section{The Deformation Theory of Semigroupal and Monoidal Functors and 
Braided Monoidal Categories}

	We now turn our attention to the actual structure of deformations of
braided monoidal categories.  As we have stated above, we proceed by the 
indirect route of developing a deformation theory for monoidal functors,
then applying it to classify deformations of a multiplication on a tensor
category. 

	As in \cite{CYdef} our starting point is the observation that the
coherence condition for the structure at hand is formally a cohomological
condition written multiplicatively. As was the case in \cite{CYdef} for
the structure maps of a deformation of a tensor category, if we consider
the structure maps for a purely functorial first order deformation of a 
semigroupal functor, we find that the main coherence condition is equivalent
(in the presence of the coherence for the map being deformed) to a
familiar cohomological condition, or rather, to what would be a familiar
condition were it not for the intervening $\lceil \;\; \rceil$'s.

	In the present case, if $\tilde{F}^\prime = \tilde{F} 
+ F^{(1)}\epsilon$ for
$F^{(1)}$ a natural transformation from $F(-\otimes -)$ to $F(-)\otimes F(-)$
and $\epsilon^2 = 0$, then the
hexagon on $\tilde{F}^\prime$ is equivalent to the condition that

\[ \lceil F^{(1)}_{A\otimes B, C} \rceil + \lceil F^{(1)}_{A,B}\otimes
F(C) \rceil = \lceil F^{(1)}_{A,B\otimes C} \rceil + \lceil F(A)\otimes 
F^{(1)}_{B,C} \rceil \]

\noindent recalling that here $\lceil \;\; \rceil$ denotes ``padding''
with coherence maps for $F$ and the semigroupal structures, as is permitted
by the results of \cite{epstein}. Of course, solving for $0$ gives us
the condition that $\delta(F^{(1)}) = 0$, where $\delta$ is the familiar
coboundary for the bar resolution, but dressed up with $\lceil \;\; \rceil$'s
on its terms, and using $\otimes $ in place of multiplication.

We now wish to describe a setting in which this ``coboundary'' is in fact the
coboundary in a cochain complex associated to the semigroupal functor $F$.
As in \cite{CYdef} let $^n \otimes$ (resp. $\otimes ^n$) denote the completely
left (resp. right) parenthesized $n$-fold iterated semigroupal product.
It should then be observed that all of the terms of the formula above
are natural transformations from $F(^3 \otimes)$ to $\otimes ^3(F^3)$.
It is thus reasonable to make

\begin{defin}
The {\em deformation complex} of a semigroupal functor targetted at a 
$K$-linear
tensor category $\cal C$ is the $K$-cochain
complex $(X^\bullet(F),\delta)$ for which 

\[ X^n(F) = Nat(F(^n \otimes), \otimes^n (F^n)) \]

\noindent with coboundary $\delta$ given by the usual bar-resolution 
formula with $\otimes$ as multiplication, and $\lceil \;\; \rceil$ 
applied to all terms. Denote the cohomology of the complex by $H^\bullet(F)$
\end{defin}

It follows from the usual argument (and the bilinearity of the composition of
natural transformations and $\otimes$ and the coherence conditions for
$F$ and $\otimes$) that $\delta^2 = 0$.
It is also easy to see that a semigroupal natural isomorphism which reduces
to the identity mod $\epsilon$ between two
first order deformations is precisely a 1-cochain which cobounds the difference
between the two 2-cocycles which name the deformations.  Thus we have

\begin{thm} \label{h2classifies}
There is a natural 1-1 correspondence between the first order 
purely functorial deformations
of a semigroupal functor $F$ and the 2-cocycles of the deformation complex of
$F$.  Moreover, the semigroupal natural isomorphism classes of first order
purely functorial
deformations of $F$ are in natural 1-1 correspondence with $H^2(F)$.
\end{thm}

Thus, $F$ admits non-trivial purely functorial deformations precisely when
$H^2(F) \neq 0$.  Because this rigidity is limited to purely functorial
deformations, we will say a semigroupal functor 
$F$ is {\em purely rigid} whenever
$H^2(F) = 0$.

We can now consider the question of extending an $n^{th}$ order deformation
to an $(n+1)^{st}$ order deformation, that is, suppose we have an $n^{th}$
order deformation

\[ F_n = F + F^{(1)}\epsilon + . . . + F^{(n)}\epsilon^n \]

\noindent where $\epsilon^{n+1} = 0$, and we wish to find an $n+1^{st}$ order
deformation

\[ F_n = F + F^{(1)}\epsilon + . . . + F^{(n+1)}\epsilon^{n+1} \]

\noindent (for $\epsilon^{n+2} = 0 \neq \epsilon^{n+1}$) by adding another term.

It is easy to see that the required condition on $F^{(n+1)}$ is that

\[ \delta(F^{(n+1)}) = \sum_{l+m = n+1, 1 \leq l,m \leq n}
\lceil F^{(m)}_{A,B}\otimes F(C) (F^{(l)}_{A\otimes B,C}) \rceil
- \sum_{p+q = n+1, 1 \leq p,q \leq n} \lceil F(A)\otimes F^{(q)}_{B,C}
(F^{(p)}_{A,B\otimes C}) \rceil \]

Thus this right-hand side may be viewed as the obstruction to extending
the $n^{th}$ order deformation to an $(n+1)^{st}$ order deformation.
As in \cite{G2}, we will see that this obstruction is always a 3-cocycle,
and thus the obstruction to extending a purely functorial deformation
by one order may be viewed as a cohomology class in $H^3(F)$.

In fact, the cochain complex associated to a semigroupal functor shares
many of the properties of the Hochschild complex of an associative 
algebra $A$ with coefficients in $A$ which were described by Gerstenhaber
\cite{G1,G2,GS}. In particular, we have two products defined on cochains:

\[ -\cup - : X^n(F) \times X^m(F)\rightarrow X^{n+m}(F) \]

\noindent given by

\[ G\cup H_{A_1, \ldots A_{n+m}} = \lceil 
F(A_1)\otimes \ldots \otimes F(A_n) \otimes
H_{A_{n+1}, \ldots ,A_{n+m}}(G_{A_1, \ldots A_n}\otimes F(A_{n+1})\otimes 
\ldots \otimes F(A_{n+m})) \rceil \]

\noindent and 

\[\langle -,- \rangle
:X^n(F)\times X^m(F)\rightarrow X^{n+m-1}(F) \]

\noindent given by

\begin{eqnarray*} \lefteqn{\langle G,H \rangle_{A_1, \ldots A_{n+m-1}} 
= } \\
 & & \sum (-1)^{mi} \lceil F(A_1)\otimes \ldots F(A_i) \otimes H_{A_{i+1},
\ldots ,A_{i+n}}\otimes F(A_{i+n+1})\otimes \ldots \otimes F(A_{n+m-1}) \\
 & & \hspace*{2cm}
(G_{A_1,\ldots ,A_i,A_{i+1}\otimes \ldots \otimes A_{i+n}, A_{i+n+1},\ldots
A_{n+m-1}}) \rceil 
\end{eqnarray*}

\begin{propo}
This latter product comes from a ``pre-Lie system''  in the terminology
of Gerstenhaber \cite{G1} given by

\begin{eqnarray*}
\lefteqn{\langle G,H \rangle^{(i)}_{A_1, \ldots A_{n+m-1}} 
= } \\
 & &  \lceil F(A_1)\otimes \ldots F(A_i) \otimes H_{A_{i+1},
\ldots ,A_{i+n}}\otimes F(A_{i+n+1})\otimes \ldots \otimes F(A_{n+m-1}) \\
 & & \hspace*{2cm}
(G_{A_1,\ldots ,A_i,A_{i+1}\otimes \ldots \otimes A_{i+n}, A_{i+n+1},\ldots
A_{n+m-1}}) \rceil 
\end{eqnarray*}

\noindent where $X^n(F)$ has degree $n-1$.
\end{propo}

\noindent{\bf proof:} 
First, note that the ambiguities of parenthesization in the semigroupal 
products
in this definition are rendered irrelevant by the $\lceil \;\; \rceil$ on
each term by virtue of the coherence theorem
for semigroupal functors.

It is obvious that the product is given by a sum of
these terms, so the content of the proposition is really that the
$\langle - , - \rangle^{(i)}$'s satisfy the definition of a pre-Lie system, 
that is for $G \in X^m(F)$, $H \in X^n(F)$, and $K \in X^p(F)$ that

\[ \langle \langle G,H \rangle^{(i)}, K \rangle^{(j)} = \left\{
	\begin{array}{ll}
	 \langle \langle G,K \rangle^{(j)}, H \rangle^{(i+p-1)} & \mbox{if 
		$0 \leq j \leq i-1$} \\
	  \langle G, \langle H, K \rangle^{(j-i)} \rangle^{(i)} & \mbox{if
		$i \leq j \leq n$}
	\end{array} \right. \]

\noindent (recall that a $k$-chain has degree $k-1$).

This is a simple computational check.  One must remember that naturality
will allow one to commute the prolongations of $K$ and $H$ in verifying
the first case.
\smallskip

We can then write the obstruction equation as 

\[ \delta(F^{(n+1)}) = \sum_{i=1}^n \langle F^{(i)}, F^{n-i+1} \rangle \]

In particular, we have a 

\begin{propo} \label{obstructions.closed}
For any $n^{th}$ order purely functorial deformation of a semigroupal
functor $F$ targetted at a $K$-linear tensor category, the obstruction
to extension to an $(n+1)^{st}$ order deformation is a 3-cocycle.
\end{propo}

\noindent{\bf proof:} This follows immediately from the proof in
Gerstenhaber \cite{G2} of \S 5, Proposition 3, 
since that proof depends only on the fact that
the operation comes from a pre-Lie system.
\smallskip

In particular, if $H^3(F) = 0$, any first order purely functorial
deformation of $F$ can
be extended to an $n^{th}$ order deformation for any $n$, and thus to a
formal series deformation.

The case of a monoidal functor is somewhat less clean than for semigroupal
functors.  Here we will only consider the case of those purely functorial
deformations for which
the unit isomorphism $F_0$ is deformed trivially, that is,
is given by
$F_0 \otimes 1$, where
$1$ is the unit of the ring $R$ over which the deformation takes place.
We call such deformations {\em proper}.

It appears to be possible to remove this restriction in the case of 
monoidal functors which are faithful and monic on objects (or equivalent
to such a functor), but we will not pursue this further in this paper.

A brief consideration of the unit coherence conditions for monoidal
functors shows that a proper $n^{th}$ order deformation of a monoidal
functor must be given by a semigroupal deformation 

\[ \tilde{F}^\prime = \tilde{F}^\prime + F^{(1)}\epsilon + \ldots
	+ F^{(n)}\epsilon^n \]

\noindent for which $F^{(k)}_{A,I} = 0$ and $F^{(k)}_{I,A} = 0$ for
all $k$ and all objects $A$.
 
This leads to the following 

\begin{defin}
The {\em proper deformation complex} of a monoidal functor targetted at
a $K$-linear tensor category is the $K$-cochain complex $(C^\bullet(F),\delta)$
for which

\[ C^n(F) = \{\phi_{A_1, \ldots ,A_n} | A_i = I \; \vdash \; 
\phi_{A_1, \ldots ,A_n} = 0\} \subset Nat(F(^n\otimes), \otimes^n(F^n)) \]

\noindent with coboundary $\delta$ given by the restriction of the coboundary
on the deformation complex of $F$.  We denote the cohomology of this 
complex by ${\bf H}^\bullet(F)$.
\end{defin}

The following proposition ensures that this definition makes sense and that the
proper deformation complex has the same type of
structures as the deformation complex:

\begin{propo}
If $G$ and $H$ are proper cochains, then so are $\delta(G)$, $G\cup H$,
and $\langle G, H \rangle^{(i)}$.
\end{propo}

\noindent{\bf proof:} For the latter two, observe that any $I$ will 
give a factor in the expression for the product in which $I$ occurs
either as an index of $G$ or as an index of $H$, and thus the product 
will be $0$.  For the first, note that all terms of $\delta(G)$ except
two will be of the form $\lceil G \rceil$ for an instance of $G$ with
$I$ as one of the indices.  The remaining two terms will be of the
same form (though without $I$ indices), of opposite signs,
and will differ only in which indices
are tensored by $I$. Thus by the presence of the $\lceil \; \; \rceil$
they will be the same and will cancel.
\smallskip

Notice in particular that the obstructions to extending a proper
deformation will always be proper cochains.

It then follows immediately that we have results corresponding to
Theorem \ref{h2classifies} and Proposition \ref{obstructions.closed}.

\begin{thm}
There is a natural 1-1 correspondence between the first order 
proper deformations
of a monoidal functor $F$ and the 2-cocycles of the proper deformation 
complex of
$F$.  Moreover, the monoidal natural isomorphism classes of first order
proper
deformations of $F$ are in natural 1-1 correspondence with ${\bf H}^2(F)$.
\end{thm}

\begin{propo}
For any $n^{th}$ order proper deformation of a monoidal
functor $F$ targetted at a $K$-linear tensor category, the obstruction
to extension to an $n+1^{st}$ order deformation is a 3-cocycle.
\end{propo}

We may now speak of a monoidal functor as {\em properly rigid} if 
${\bf H}^2(F) = 0$, and observe that if ${\bf H}^3(F) = 0$ then any
first order proper deformation can be extended to a formal series
deformation.

Now, consider the case of a multiplication on a tensor category $\cal C$.
The requirement that there exist right and left unit 
transformations for the multiplication imposes an additional condition on
the admissible deformations.  However, these also are readily understood
in terms of the cohomology theory. 

 In particular, notice that if $F$ and
$G$ are a composable pair of monoidal (or semigroupal) functors, 
a deformation of either
leads to an induced deformation of the composite (provided the appropriate
target categories are $K$-linear tensor categories).  More precisely,
if $\tilde{F}$ is a deformation of $F$, then $\tilde{F}(G)$ is a 
deformation of $F(G)$, while if $\tilde{G}$ is a deformation of $G$, and
both the source and target of $F$ are $K$-linear tensor categories, then
$\tilde{F}_{\rm trivial}(\tilde{G}_!)$ is a deformation of $F(G)$, 
where $\tilde{F}_{\rm trivial}$ denotes the trivial deformation of $F$
and $\tilde{G}_!$ denotes the $R$-linear extension of $\tilde{G}:{\cal C}
\rightarrow \tilde{\cal D}$ to $\tilde{\cal C}$. In terms
of induced deformations we then have

\begin{thm}
If $\Phi$ is a multiplication on a tensor category $\cal C$, and 
$\Phi^\prime$ is a proper
deformation of $\Phi$ as a monoidal functor
then  $\Phi^\prime$
is a multiplication on $\tilde{\cal C}$ precisely when the induced 
deformations $\Phi^\prime (-,I)$ and $\Phi^\prime (I,-)$
are trivial.
\end{thm}

\noindent{\bf proof:} If the induced deformations are trivial, we can
compose the trivializing natural isomorphsism with the $R$-bilinear extension
of $\frak r$ (resp. $\frak l$) to obtain the deformed $\frak r$ (resp. $\frak
l$).  Conversely, if there are deformed $\frak r$ and $\frak l$, their
composition with the $R$-bilinear extension of the inverse to the undeformed
$\frak r$ and $\frak l$ give the respective trivializations.
\smallskip

\section{Vassiliev Invariants for Framed Links}

	It is usual to discuss Vassiliev theory in terms of unframed 
oriented knots and links (cf. \cite{BL}, \cite{Stanford}, \cite{Vas}).  
We will, however, remain in the setting 
most natural for functorial invariants (and incidentally most closely
connected to 3- and 4-manifold topology), that of framed links.

	Following Goryunov \cite{Gor} one can regard framed knots (and
links) not as ordinary knots (and links) with additional structure, but
as equivalence classes of mappings from open annular neighborhoods $U$ of
$S^1$ (or disjoint unions of such) into ${\Bbb R}^3$.  Observe that
any such mapping $g:U\rightarrow {\Bbb R}^3$ induces a mapping $Tg$
from $i^\ast(T{\Bbb R}^2)$ to
$T{\Bbb R}^3$. (Here $i$ is the inclusion of $S^1$ or a disjoint union
of $S^1$'s into $U$.)

	For ease in the link setting, we need to specify a countable 
family of disjoint circles:  say 
$S_n = \{ z \; | \; |z-3n| = 1\} \subset {\Bbb C}$, then consider
disjoint annular neighborhoods of $S_1 \cup S_2 \cup ... \cup S_k$
when dealing with $k$-component links.

\begin{defin} Two mappings $g_i:U_i\rightarrow {\Bbb R}^3$ $i=1,2$ are
{\em equivalent} if the $U_i$'s are annular neighborhoods of the
same $S_1 \cup ... \cup S_k$ and
the mappings $Tg_i:i^\ast(T{\Bbb R}^2)\rightarrow T{\Bbb R}^3$
coincide.
\end{defin}

The space of (possibly singular) framed links is the space of all
$C^\infty$ mappings of $U$'s modulo this equivalence relation. We denote
it by $\Omega_f$.  The subspace of equivalence classes of mapping of 
neighborhoods of $S_1 \cup ... \cup S_k$ is the space of
(possibly singular) framed links of $k$ components, and will be
denoted $\Omega_f(k)$. Now consider the subspace ${\cal O}_f$
of all (equivalence
classes of) mappings such that $Tg$ is an embedding.

\begin{defin}
A {\em (non-singular) framed link} is a connected component of ${\cal O}_f$.
A connected component of ${\cal O}_f(k) = {\cal O}_f \cap \Omega_f(k)$ is
a {\em (non-singular) framed link of $k$ components}.
\end{defin}

First, note that we will often drop the adjective ``non-singular'' to match the
usual usage in knot theory (as in earlier sections).
Second, the designation of these maps as non-singular implicitly identifies
the element of $\Omega_f \setminus {\cal O}_f$ as {\em singular}.
We denote the discriminant locus $\Omega_f \setminus {\cal O}_f$ by
$\Sigma_f$. As observed in Goryunov \cite{Gor} the discriminant is the
union of two hypersurfaces, one on which the framing degenerates, denoted
$\Sigma_f^\prime$, and one on which the disjoint union of circles is not
embedded, denoted $\Sigma_f^{\prime \prime}$. 

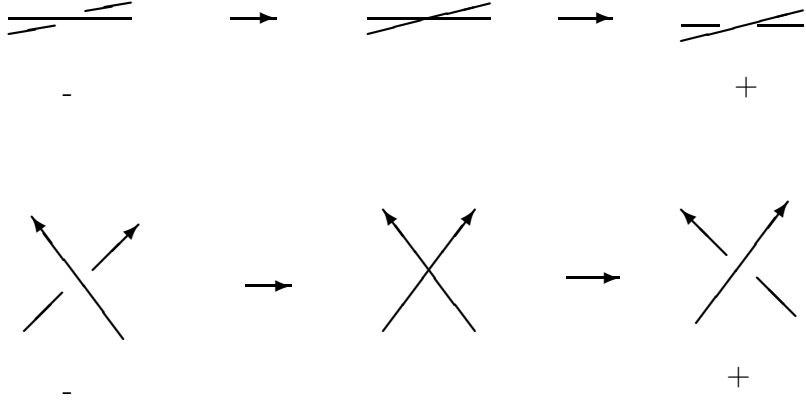
\begin{figure}[hbt]
\setlength{\unitlength}{0.008in}%
\begin{center} \begin{picture}(520,260)(90,520)
\thicklines
\put(325,770){\line( 1, 0){ 80}}
\put(325,760){\line( 4, 1){ 80}}
\put(540,570){\vector( 3, 4){ 60}}
\put(560,615){\vector(-1, 1){ 30}}
\put(580,600){\line( 1,-1){ 25}}
\put(165,560){\vector(-3, 4){ 60}}
\put(145,605){\vector( 1, 1){ 30}}
\put(125,590){\line(-1,-1){ 25}}
\put( 90,770){\line( 1, 0){ 80}}
\put( 90,760){\line( 6, 1){ 30}}
\put(140,775){\line( 6, 1){ 30}}
\put(530,755){\line( 4, 1){ 80}}
\put(530,765){\line( 1, 0){ 25}}
\put(580,765){\line( 1, 0){ 30}}
\put(235,770){\vector( 1, 0){ 30}}
\put(450,770){\vector( 1, 0){ 35}}
\put(335,565){\vector( 3, 4){ 60}}
\put(395,565){\vector(-3, 4){ 60}}
\put(245,595){\vector( 1, 0){ 30}}
\put(455,600){\vector( 1, 0){ 35}}
\put(565,720){\makebox(0,0)[lb]{\raisebox{0pt}[0pt][0pt]{\twlrm +}}}
\put(125,715){\makebox(0,0)[lb]{\raisebox{0pt}[0pt][0pt]{\twlrm -}}}
\put(125,520){\makebox(0,0)[lb]{\raisebox{0pt}[0pt][0pt]{\twlrm -}}}
\put(560,530){\makebox(0,0)[lb]{\raisebox{0pt}[0pt][0pt]{\twlrm +}}}
\end{picture} \end{center}

\caption{Co\"{o}rienting the finite codimensional strata \label{coorient}}
\end{figure}

The minimal degenerations of each type are illustrated in Figure 
\ref{coorient}.  We denote the subspace of $\Sigma_f$ in which there
are exactly $n$ degenerations of either type by $\Sigma_{f,n}$, and its
intersection with $\Omega_f(k)$ by $\Sigma_{f,n}(k)$.

As in Goryunov \cite{Gor}, we co\"{o}rient the finite codimensional
strata of $\Sigma_f$ by the local prescriptions given in Figure
\ref{coorient}, and give a Vassiliev type prescription for the
extension of invariants of framed knots to singular framed knots as in
Figure \ref{extend}.

\begin{figure}[hbt]
\setlength{\unitlength}{0.008in}%
\begin{center} \begin{picture}(601,195)(45,560)
\thicklines
\put(555,560){\vector( 3, 4){ 60}}
\put(615,560){\vector(-3, 4){ 60}}
\put(100,560){\vector( 3, 4){ 60}}
\put(120,605){\vector(-1, 1){ 30}}
\put(140,590){\line( 1,-1){ 25}}
\put(380,560){\vector(-3, 4){ 60}}
\put(360,605){\vector( 1, 1){ 30}}
\put(340,590){\line(-1,-1){ 25}}
\put( 90,710){\line( 4, 1){ 80}}
\put( 90,720){\line( 1, 0){ 25}}
\put(140,720){\line( 1, 0){ 30}}
\put(315,720){\line( 1, 0){ 80}}
\put(315,710){\line( 6, 1){ 30}}
\put(365,725){\line( 6, 1){ 30}}
\put(540,720){\line( 1, 0){ 80}}
\put(540,710){\line( 4, 1){ 80}}
\put( 80,723){\oval( 20, 64)[tl]}
\put( 80,723){\oval( 20, 66)[bl]}
\put( 80,598){\oval( 20, 64)[tl]}
\put( 80,598){\oval( 20, 66)[bl]}
\put(300,723){\oval( 20, 64)[tl]}
\put(300,723){\oval( 20, 66)[bl]}
\put(530,723){\oval( 20, 64)[tl]}
\put(530,723){\oval( 20, 66)[bl]}
\put(300,598){\oval( 20, 64)[tl]}
\put(300,598){\oval( 20, 66)[bl]}
\put(535,598){\oval( 20, 64)[tl]}
\put(535,598){\oval( 20, 66)[bl]}
\put(180,723){\oval( 20, 66)[br]}
\put(180,723){\oval( 20, 64)[tr]}
\put(185,598){\oval( 20, 66)[br]}
\put(185,598){\oval( 20, 64)[tr]}
\put(410,598){\oval( 20, 66)[br]}
\put(410,598){\oval( 20, 64)[tr]}
\put(410,723){\oval( 20, 66)[br]}
\put(410,723){\oval( 20, 64)[tr]}
\put(635,723){\oval( 20, 66)[br]}
\put(635,723){\oval( 20, 64)[tr]}
\put(635,598){\oval( 20, 66)[br]}
\put(635,598){\oval( 20, 64)[tr]}
\put(230,715){\makebox(0,0)[lb]{\raisebox{0pt}[0pt][0pt]{\twlrm -}}}
\put(230,595){\makebox(0,0)[lb]{\raisebox{0pt}[0pt][0pt]{\twlrm -}}}
\put(465,715){\makebox(0,0)[lb]{\raisebox{0pt}[0pt][0pt]{\twlrm =}}}
\put(460,595){\makebox(0,0)[lb]{\raisebox{0pt}[0pt][0pt]{\twlrm =}}}
\put( 45,715){\makebox(0,0)[lb]{\raisebox{0pt}[0pt][0pt]{\twlrm V}}}
\put(275,710){\makebox(0,0)[lb]{\raisebox{0pt}[0pt][0pt]{\twlrm V}}}
\put(505,710){\makebox(0,0)[lb]{\raisebox{0pt}[0pt][0pt]{\twlrm V}}}
\put( 55,585){\makebox(0,0)[lb]{\raisebox{0pt}[0pt][0pt]{\twlrm V}}}
\put(270,595){\makebox(0,0)[lb]{\raisebox{0pt}[0pt][0pt]{\twlrm V}}}
\put(505,590){\makebox(0,0)[lb]{\raisebox{0pt}[0pt][0pt]{\twlrm V}}}
\end{picture} \end{center}

\caption{Vassiliev type extension formulae \label{extend}}
\end{figure}

\begin{defin}
A Vassiliev invariant of links (resp. links of $k$ components) is a 
locally constant function on ${\cal O}_f$ (resp. ${\cal O}_f(k)$) whose
extension according to the prescription of Figure \ref{extend} vanishes
on $\Sigma_{f,n+1}$ (resp. $\Sigma_{f,n+1}(k)$) for some $n$, in which case
the invariant
is said to be of type $\leq n$.  The invariant is of type $N$ when $N$ is the
minimal such $n$.
\end{defin}

\section{Deformations and Vassiliev Invariants}

We are now in a position to state our main theorem:

\begin{thm} \label{main.thm}
Let $\cal C$ be any $K$-linear rigid symmetric monoidal category, and let
$\tilde{\cal C}$ be any $n^{th}$ order tortile deformation of $\cal C$.
For any object $X$ of $\tilde{\cal C}$, let $V_X$ denote the functor
from $\cal FT$ to $\tilde{\cal C}$ induced by Shum's Coherence Theorem.
Then $V_X$ restricted to $End(I)$, regarded as the set of framed links,
is a $K[\epsilon]/<\epsilon^{n+1}>$-valued
Vassiliev invariant of type $\leq n$, and is, moreover, multiplicative
under disjoint union.
\end{thm}

From this will follow, almost as a corollary

\begin{thm} \label{second.thm}
Let $\cal C$ be any $K$-linear rigid symmetric monoidal category, and let
$\tilde{\cal C}$ be any $n^{th}$ order tortile deformation of $\cal C$
(resp. formal series deformation of $\cal C$).
For any object $X$ of $\tilde{\cal C}$, let $V_X$ denote the functor
from $\cal FT$ to $\tilde{\cal C}$ induced by Shum's Coherence Theorem,
and let $V_{X,k}$ denote the $K$-valued framed link invariant which assigns to
any link the coefficient of $\epsilon^k$, for $k = 0,...n$, (resp. 
for $k \in {\Bbb N}$).
Then $V_{X,k}$ is a
Vassiliev invariant of type $\leq k$.
\end{thm}

	The proof of Theorem \ref{main.thm} is quite simple and 
similar to previous proofs of similar results:

	The key is to observe that the bilinearity of composition in
$\tilde{\cal C}$ allows us to use the Vassiliev prescription to extend
the functor $V_X$ from $\cal FT$ to a larger category of
singular framed tangles, $\widetilde{\cal FT}$, 
whose maps are isotopy classes of framed tangles with finitely many
degeneracies of either of the two basic types.

	Consider a singular framed link with $n+1$ degeneracies (of either
type).  Now, we can represent the framed link as a composition of singular
framed tangles, each of which has at most one degeneracy, crossing, framing
twist, or
extremum. Now
the value of the extended functor on such a tangle with a
 degeneracy of the first type
(framing degeneracy) is a monoidal product of identity maps with
$\theta_X - \theta_X^{-1}$ (or its dual)   ,
while the value on such a tangle with a degeneracy of the second type is
a monoidal product of identity maps with $\sigma_{X,X} - \sigma_{X,X}^{-1}$.

Now, observe that
 
\[ \theta_X - \theta_X^{-1} \in Hom_{\cal C}(X,X)\otimes <\epsilon> \]

\noindent and

\[ \sigma_{X,X} - \sigma_{X,X}^{-1} \in Hom_{\cal C}(X\otimes X, X\otimes X)
\otimes <\epsilon>.\]

It follows from the bilinearity of composition in $\tilde{\cal C}$ that
the composite representing the framed link as an element of $End(I)$ lies
in $End_{\cal C}\otimes <\epsilon^{n+1}> = 0$, 
thus showing $V_X$ to be Vassiliev of type
$\leq n$. Multiplicativity follows from functoriality.

	Theorem \ref{second.thm} follows from the Theorem \ref{main.thm},
Lemma \ref{reduction}, and the following lemma, the proof of which
is a trivial exercise:

\begin{lemma}
If $V$ is an $A$-valued Vassiliev invariant, 
and $f:A\rightarrow B$ is a linear map ($A$ and $B$ here are abelian groups)
then $f(V)$ is a $B$-valued Vassiliev invariant.
\end{lemma}

	On general principles these results are more general than previous
results of this form:  our initial category can be any $K$-linear rigid
symmetric monoidal category for any field $K$, not just the category of 
representations of
a (semi)simple Lie algebra.  In particular, this means that the invariants
obtained from any representation of a $q$-deformed universal enveloping
algebra for a super-Lie algebra (cf. \cite{Man}) 
give rise to sequences of Vassiliev
invariants for links of any number of components when we set $q = e^h$ 
(or any other formal series for that matter) and
consider the coefficients of $h^n$. (The observation about using other series,
of course applies in the ordinary Lie algebra case as well.)  
Likewise the coefficients of any 
$h^n$ in invariants arising from a representation of a multiparameter quantum
group (cf. \cite{DD}, \cite{Haz}, \cite{Tak})and {\em any} instantiation 
of the deformation parameters by power series will be Vassiliev invariants
(even in the case of quasi-Hopf deformations).  In fact, by repeating 
the arguments
given above for quotients of $K[[x]]$ by nilpotent ideals in the case of
$K[[h_1,...,h_n]]$, one can show that the coefficients of any monomial
in $h_1,...,h_n$ in the invariant associated to a representation of
a multiparameter quantum group and any instantiation of the $i^{th}$ 
deformation parameter by a power series in $h_i$ are Vassiliev invariants.

\section{Questions Raised}

The results contained herein are only partial.  In order to fully relate
Vassiliev theory to the deformation theory of braided monoidal categories,
it will be necessary to give a good accounting of that deformation theory
in the non-proper case.  The difficulties here are two-fold.  First is
the removal of the restriction that the unit isomorphism not be deformed,
and second is the removal of the purely functorial restriction.

We thus ask:

\begin{enumerate}
\item Is there a good cohomological deformation theory for purely
functorial (but not necessarily proper) deformations of monoidal functors?
\item Can the deformation theory of this work be combined in a satisfactory
way with the deformation theory of monoidal categories given in
\cite{CYdef} to provide a deformation theory for (not necessarily purely
functorial) deformations of monoidal functors?
\item Even a satisfactory answer to question 2.\ will not suffice to give
a complete deformation theory for braided monoidal categories, since in
2.\ only the target is deformed, while a deformationthe structure maps of the 
underlying category, induce deformations of
the structure maps of both the source
and the target of the multiplication.  Thus, we ask:
Is there a satisfactory cohomological deformation theory for unrestriced
deformations of multiplications on monoidal categories?
\end{enumerate}

Likewise, the construction given here of Vassiliev invariants from 
deformations of monoidal categories leads to a series of questions 
regarding the
converse problem: which Vassiliev invariants arise from functorial invariants?

Chief among these are various precise formulations of that basic question:

\begin{enumerate}
\item Given a $K$-valued
Vassiliev invariant $v_n$ of type $n$ of framed knots, when
do there exist Vassiliev invariants $v_0,...,v_{n-1}$ of types $0,...,n-1$
respectively, such that the $K[\eta]/<\eta^n+1>$-valued invariant

\[ v_0 + v_1 \eta + ...+ v_n \eta^n \]

\noindent is the restriction to framed knots
of an invariant defined by an object in
some tortile deformation of a $K$-linear rigid symmetric monoidal category?
\item  Given a $K$-valued Vassiliev invariant $v$ of
type $n$ of $m$-component links,
\begin{enumerate}
\item When does $v$ extend to a Vassiliev invariant of links as defined
above?
\item When do there exist Vassiliev invariants of lower orders so that
the generating function as in question 1.\ is the invariant defined by an object
in some tortile deformation of a $K$-linear rigid symmetric monoidal
category?
\end{enumerate}
\item Given a Vassiliev invariant of links as described above, when can
it be completed as in questions
1.\ and 2.\ by adding lower order Vassiliev invariants
to give a functorial invariant?

\end{enumerate}

	A partial result toward answers to these three questions
can be read off from the results
above, namely, a necessary condition for a Vassiliev invariant of framed
links to arise as a series coefficient in a tortile deformation of a 
symmetric tensor category:

\begin{thm} \label{disjoint.restriction}
If $v$ is a Vassiliev invariant of (framed) links arising
as the coefficient of $\epsilon^n$ in an $n^{th}$ order tortile deformation
of a rigid symmetric tensor category, then there exist Vassiliev invariants
$v_k$ of (framed) links $k = 1, . . . ,n$ (with $v_n =v$)  such that 

\[ v(K \coprod L) = v_0(K)v_n(L) + v_1(K)v_{n-1}(L) + . . . v_n(K)v_0(L) \]

\noindent and each $v_i$ is of type $\leq i+1$, where $\coprod$ denotes
separated union of links.
\end{thm}

	In terms of Vassiliev invariants for links of a fixed number of
components, we can similarly see

\begin{thm}
If $v$ is a Vassiliev invariant of (framed) links of $N$-components
arising as the coefficient of $\epsilon^n$ in a tortile deformation of a rigid
symmetric tensor category, then for all $M = 1, . . . ,N$ and $i = 0, . . . ,n$
there exist Vassiliev invariants $v_{M,i}$ of $M$ component (framed) 
links of type $\leq i$
such that for all $K \coprod L$ separated unions of links, with
$K$ having $M$ components and $L$ having $N-M$ components, we have

\[ v(K \coprod L) = v_{M,0}(K)v_{N-M,n}(L) + v_{M,1}(K)v_{N-M,n-1}(L) + . . .
	+ v_{M,n}(K)v_{N-M,0}(L) \]
\end{thm}
 
	Finally, there are questions of deeper relationships between the
cohomological deformation theories of
 of monoidal categories and of (bi)algebras and Vassiliev theory:

\begin{enumerate}
\item What is the relationship between the deformation
theory of braided monoidal categories and the original 
cohomological formulation of
Vassiliev invariants \cite{Vas}?   
\item What relationship is there between the Vassiliev invariants
arising from the tortile deformations of a particular symmetric monoidal
category and the structure of the category itself?
\item What is the relationship between the deformation theory of a (bi)algebra
$A$
and that of monoidal functors associated to the algebra (in particular, the
forgetful functor 

\[ U:Rep(A)\rightarrow K-vect, \] 

\noindent and for (quasi-)triangular
$A$ the tensor product as a monoidal functor with stuctural transformations
as induced by the braiding or symmetry as in Theorem \ref{mult.from.braiding}?

\end{enumerate}

\end{document}